\documentclass[11pt,leqno,a4paper]{article}
\usepackage{amsmath}
\usepackage{amssymb}
\advance\topmargin by -2.2cm
\newcommand{\bbr}{I\!\!R}

\newcommand{\bbn}{I\!\!N}

\newcommand{\cala}{{\cal A}}
\newcommand{\calb}{{\cal B}}
\newcommand{\calc}{{\cal C}}

\newcommand{\calf}{{\cal F}}

\newcommand{\call}{{\cal L}}

\newcommand{\barr}{\begin{array}}
\newcommand{\earr}{\end{array}}
\newcommand{\beqq}{\begin{equation}}
\newcommand{\eeqq}{\end{equation}}
\newcommand{\beao}{\begin{eqnarray*}}
\newcommand{\eeao}{\end{eqnarray*}\noindent}
\newcommand{\beam}{\begin{eqnarray}}
\newcommand{\eeam}{\end{eqnarray}\noindent}

\newcommand{\halmos}{\quad\hfill\mbox{$\Box$}}

\newcommand{\la}{\lambda}

\newcommand{\si}{\sigma}
\newcommand{\al}{\alpha}

\newcommand{\vp}{\varphi}

\newcommand{\ov}{\overline}
\newcommand{\wh}{\widehat}
\newcommand{\wt}{\widetilde}

\newcommand{\lra}{\longrightarrow}

\newcommand{\nto}{n\to\infty}

\newcommand{\var}{{\rm Var}}

\begin{document}

{\Large\bf A time inhomogenous Cox-Ingersoll-Ross diffusion with jumps}   \\

Reinhard H\"opfner, Universit\"at Mainz\\

\vskip1.5cm\small
{\bf Abstract: } 
We consider a time inhomogeneous Cox-Ingersoll-Ross diffusion with positive jumps. We exploit a branching property to prove existence of a unique strong solution under a restrictive condition on the jump measure. We give Laplace transforms for the transition probabilities, with an interpretation in terms of limits of mixtures over Gamma laws. 
\normalsize

\vskip1.2cm
We look at a stochastic differential equation of Cox-Ingersoll-Ross (CIR) type with jumps
$$
d\xi_t \;=\; \left[ a(t) - \beta(t)\, \xi_{t^-}\right] dt  \;+\; \int_{(0,\infty)} y\; \mu(dt,dy) 
\;+\; \si(t)\,\sqrt{\xi_{t^-}\!\!\vee 0}\; dW_t  \;,\quad t\ge 0   
\leqno(\diamond)
$$  
which is inhomogeneous in time.  
Some starting point $\xi_0\equiv x_0\ge 0$ is fixed. We are interested in strong solutions driven by $(W,\mu)$ where  Brownian motion $W$ and Poisson random measure (PRM) $\mu$ are assumed independent, such that
$$
\mu(ds,dy) \;\;\mbox{on}\;\; (0,\infty){\times}(0,\infty) \;\;\mbox{has intensity}\;\; \wt a(s) ds\, \nu(dy)  
$$
for some $\si$--finite measure $\nu$ on $(0,\infty)$ with  
\beqq\label{permanentcondition}
\int_0^\infty (y\wedge 1)\, \nu(dy) \;<\;  \infty \;. 
\eeqq
Hence  jumps of $\xi$ are positive and summable.  The functions $a(\cdot) , \wt a(\cdot) \ge 0$ correspond to differently structured  input into equation $(\diamond)$, and there is a time-dependent backdriving force $\beta(\cdot)\ge 0$.\\

The aim of this note is to show that equation $(\diamond)$ has a unique strong solution at least for a restricted class of measures $\nu$ (theorem 1.1), and to give explicit expressions for the Laplace transforms of the transition probabilities of the process $(\diamond)$ (theorem 1.2 and remark 1.3). It will be seen that  input in terms of $a(\cdot)$ and input in terms of $\wt a(\cdot)$ induces terms of very similiar structure in the transition probabilities. Thus we extend known results  (see Ikeda and Watanabe  [IW 89, p.\ 235]) for the classical CIR diffusion 
\beqq\label{classicalCIR}
d\zeta_t \;=\; \left[ a - \beta\, \zeta_t\right] dt   
\;+\; \si\,\sqrt{\zeta_t\vee 0}\; dW_t  \;,\quad t\ge 0   
\eeqq  
(Cox, Ingersoll and Ross [CIR 85] where $a,\beta,\si >0$  are constants) to a time inhomogeneous setting with jumps. For the classical CIR diffusion (\ref{classicalCIR}), it is known that transition probabilities can be represented as Poisson mixtures of Gamma laws (combine [IW 89, p.\ 235] with Barra [B~71, p.\ 82], or see Chaleyat-Maurel and Genon-Catalot [CG 06, appendix A3]). The structure of the Laplace transform for transition probabilities in $(\diamond)$ according to theorem 1.2 and remark 1.3 below indicates additional 'mixtures' (more exactly: weak limits of additional mixtures), in time (due to time-dependence of $t\to a(t)$ and $t\to \wt a(t)$) and in space (in the sense of Poisson random measures on $(0,\infty)$ acting with intensity proportional to $\nu(dy)$ at times $t\ge 0$).  \\

One motivation for this note is an interest in statistical problems for non time homogeneous processes. In a recent work with Kutoyants [HK 09], we have a $T$-periodic deterministic function of time $t\to S(t)$  in the drift of an SDE,  all other coefficients of the SDE being time homogeneous,  and need as principal assumption --~in order to formulate limit theorems~-- a Harris condition which in order to be checked requires explicit knowledge of the transition probabilities corresponding to time steps which are integer multiples of $T$. \\

Another motivation for equation $(\diamond)$ is to model the membrane potential in a neuron which --~being part of an active network of neurons~-- receives from its network both 'continuous' and 'discrete' input. Depending on some external stimulus and hence on a  time-varying degree of activity of the network as a whole, this input is not a function of the membrane potential in the receiving neuron. In equation $(\diamond)$, we think of the 'continuous input' part $a(s)ds$ as synaptic input collected through the dendrites of the neuron, arriving decayed and delayed at the soma (cf.~[BH 06]), whereas the 'discrete input' part $\int_{(0,\infty)} y\, \mu(dt,dy)$ models synapses located  at the soma itself, or in an immediate vicinity (thus near to the recording electrode). There are data sets recording the membrane potential in one neuron (belonging to a cortical slice observed in vitro which is stimulated by a  potassium bath, cf.\  [H 07, section 3.2] or [J 09, section 5.3]; these data have been recorded by  Dr.\ W.\ Kilb, Institute of Physiology, University of Mainz) which do support a CIR type modelization of the membrane potential. Note that the classical CIR diffusion (\ref{classicalCIR}) has been used as a model for the membrane potential in many papers since L\'ansk\'y and L\'ansk\'a [LL 87]. However, only under the effect of quite  particular experimental conditions (such as electronic stabilisation of the membrane potential in the neuron which was the case for the data investigated in [J 09, section 5.3]) real data exhibit a behaviour which can safely be considered as 'stationary'. Inhomogeneity in time (one example in [H 07, section 4.6]) seems to be the general situation, and only few work has been done so far in this direction. \\

\section*{1. Results}

We state our assumptions on the coefficients of SDE $(\diamond)$ in detail. We assume the backdriving force $\beta(\cdot)$  continuous and nonnegative. In the volatility $(t,y)\to\si(t)\sqrt{y^+\,}$, the function $\si(\cdot)$ is continuous and strictly positive on $[0,\infty)$. The functions $a(\cdot)$, $\wt a(\cdot)$ modelling time dependence of the input are nonnegative, right continuous and locally bounded; they satisfy $\int_N^\infty (a+\wt a)(t)\, dt = \infty$ for arbitrarily large~$N$. $\;\wh\mu(ds,dy)$ denotes the compensator $\wt a(s)\, ds\, \nu(dy)$ of $\mu(ds,dy)$, and 
$$
\wt\mu(ds,dy) \;:=\; \mu(ds,dy) - \wh\mu(ds,dy)  
$$
the compensated random measure. Recall that $\mu$ and $W$ are independent. Our permanent assumption on $\nu$ is (\ref{permanentcondition}). Whenever $\nu$ has infinite total mass (hence $\xi$ an infinite number of small jumps over finite time intervals) our construction (strong solutions to $(\diamond)$ through an approximation which relies on a 'branching property') will require the restrictive condition
\beqq\label{restrictivecondition}  
\int_0^\infty \left(\sqrt{y}\wedge 1\right)\, \nu(dy) \;<\; \infty
\eeqq
for existence of a strong solution to equation  $(\diamond)$ driven by $W$ and $\mu$. Below, at the stages of the proofs where (\ref{restrictivecondition}) comes in, we mention this condition explicitely, working through all other parts of the proofs under the permanent assumption (\ref{permanentcondition}) which ensures summability of small jumps. \\

{\bf 1.1 Theorem: }  Condition (\ref{restrictivecondition}) implies existence of a unique strong solution for equation $(\diamond)$. \\

In equation $(\diamond)$ and in similiar equations with subintervals $I$ of $(0,\infty)$ 
$$
d\wt\xi_t \;=\; \left[ a(t) - \beta(t)\, \wt\xi_{t^-}\right] dt  \;+\; \int_I y\; \mu(dt,dy) 
\;+\; \si(t)\,\sqrt{\wt\xi_{t^-}\!\!\vee 0}\; dW_t  \;,\quad t\ge 0   
$$
which we shall consider, the process of jumps $\Delta\wt\xi = \left(\wt\xi_t-\wt\xi_{t-}\right)_{t\ge 0}$ does not depend on $\left(\wt\xi_{t-}\right)_{t\ge 0}$. Hence, under assumption (\ref{permanentcondition}) alone, 
an extension of the Yamada-Watanabe criterion ([YW 71]) on the lines of [KS 91, p.\ 291] guarantees pathwise uniqueness. An extension to more general SDE with summable jumps is given in [H 08]. Thanks to pathwise uniqueness, in exact analogy to the well known result for classical SDE due to Ikeda and Watanabe ([IW 89, theorem IV.1.1]), it will be sufficient to construct some weak solution to $(\diamond)$, or to similiar equations involving subintervals $I$ of $(0,\infty)$ as above. We construct weak solutions using a 'branching property' inherent to equation $(\diamond)$ (see lemma 2.2). It is  this 'branching property' which leads to the restrictive condition (\ref{restrictivecondition}) when --~for measures $\nu$ having infinite total mass~-- we pass from subintervals $I_n=[\delta_n,1]$ with $\delta_n\downarrow 0\,$ to $\,I=(0,1]$ and take limits. This will be the topic of section~2; we will finish the proof of theorem 1.1 in 2.6 below. \\

We turn to Laplace transforms of transition probabilities for the process $(\diamond)$. For $0\le s<t<\infty$ and $\la\ge 0$, introduce  
\beqq\label{3}
C(s,t) \;:=\; \int_s^t \frac{\si^2(v)}{2}\, \exp\left( \int_0^v \beta(u)du \right)\, dv \; , \quad
B(s,t) \;:=\; \exp\left( -\int_s^t \beta(v)dv \right) \; , 
\eeqq 
\beqq\label{3bis}
p(s,t) \;:=\; \frac{1}{B(0,t)\, C(s,t)} \; , \quad \gamma(s,t) \;:=\; \frac{1}{B(0,s)\, C(s,t)} \; , \\[2mm]
\eeqq 
\beqq\label{3ter}
\Psi_{s,t}(\la) \;:=\; \gamma(s,t)\, \frac{\la}{p(s,t)+\la} \; , \quad
\wt\Psi _{s,t}(\la) \;:=\; \int_{(0,\infty)} [ 1 - e^{-y\Psi_{s,t}(\la)}]\, \nu(dy)   \;. \\[2mm]
\eeqq
Note that $\wt\Psi _{s,t}(\la)$ is well defined under assumption (\ref{permanentcondition}) alone.\\

{\bf 1.2 Theorem: } Assume  (\ref{restrictivecondition}) in order to have a unique strong solution to $(\diamond)$. Assume  $\beta(\cdot)$ strictly positive on $[0,\infty)$. Then for $0\le s<t<\infty$ and $\la\ge 0$, 
\beao
E \left( e^{ -\, \la\, \xi_t } \mid \xi_s=y \right) \;\; = \;\;   
\exp\left( -\; y\,\Psi_{s,t}(\la) \;-\; \int_s^t  \left[\, a(v)\, \Psi_{v,t}(\la) + \wt a (v)\, \wt\Psi_{v,t}(\la) \,\right] dv \,\right) 
\eeao
is the Laplace transform of the transition probability $P\left(\xi_t\in\cdot \mid \xi_s=y\right)$ 
in the process $(\diamond)$.\\

{\bf 1.3 Remark: } In the special case $a(\cdot) = \al\, \frac{\si^2(\cdot)}{2}$ for some constant $\al>0$, we have 
$$  
\exp\left( - \int_s^t a(v)\, \Psi_{v,t}(\la)\, dv\, \right) 
\;=\; \left( 1 + \frac{\la}{p(s,t)}\right) ^{-\al} 
\;,\quad \la\ge 0  
$$ 
which is the Laplace transform of the Gamma law $\Gamma(\al,p(s,t))$. Hence, in reduction to the special case of the classical CIR diffusion (\ref{classicalCIR}), the Laplace transform given in theorem 1.2 coincides with the well known expression given in Ikeda and Watanabe [IW 89, p.\ 235]. \\
 
See also Overbeck and Ryden [OR 97]. Laplace transforms for transition probabilities in  $(\diamond)$ and in related equations  will be the topic of section 3. The interpretation of  transition probabilities in $(\diamond)$ in terms of limits of Poisson mixtures of Gamma laws will become clear through lemmata 3.3, 3.7, 3.8, 3.10. Remark 1.3 will be proved in 3.4+3.5. We will finish the proof of theorem 1.2 in 3.11 below. \\

\section*{2. Proof of theorem 1.1} 

In this section, we prove theorem 1.1. First, a localization argument shows that it is sufficient to prove theorem 1.1 in case where in addition to the assumptions stated above  
\beqq\label{localization1}
\begin{array}{l}
\mbox{$\si(\cdot)$ is bounded away from both $0$ and $\infty$ on $[0,\infty)$}\\
\mbox{$\beta(\cdot)$, $\al(\cdot)$, $\wt\al(\cdot)$ are bounded on $[0,\infty)$} \;. 
\end{array}
\eeqq
We shall assume (\ref{localization1}) throughout this proof section. Second, since pathwise uniqueness holds for equation $(\diamond)$ by the Yamada-Watanabe criterion, as in [H 08], we have to prove existence of a solution. Our key step is the following lemma. Recall that under the assumptions of theorem 1.1, the backdriving force $\beta(\cdot)\ge 0$ in $(\diamond)$ may vanish. \\

{\bf 2.1 Lemma: }  For $0\le s<\infty$ and $u>0$ deterministic, consider the equation 
$$
d\xi_t \;=\; -\beta(t)\, \xi_t\, dt  \;+\; \si(t)\,\sqrt{\xi_t^+}\; dW_t  \;, \quad t\ge s \;,\; \xi_s=u 
\leqno(*,s,u)
$$
up to time $\tau := \inf\{ t>s : \xi_t=0\}$, and define $\xi\equiv 0$ on $[[\tau,\infty[[$. Write 
$$
\xi^{(s,u)} \;=:\; \xi    
$$
for the unique strong solution. Then $\tau$ is finite a.s., and we have for fixed $T\in(0,\infty)$ 
$$
E\left( \sup_{s\le t\le s+T} (\xi^{(s,u)}_t)^2 \right)  \;\;\le\;\;   C(1+T)\; \left(\, u + u^2 \,\right)    
$$
where $C$ depends on upper and lower  bounds for $\si^2(\cdot)$ according to  (\ref{localization1}), and not on $s$, $u$, $T$. \\

{\bf Proof: } 1)~  The solution of $(*,s,u)$ is unique and strong in virtue of the Yamada-Watanabe criterion (see 
Yamada and Watanabe [YW 71], Ikeda and Watanabe [IW 89, theorem IV.1.1], Karatzas and Shreve [KS 91, p.\ 291]; note that  on every $[[0,\rho_\ell]]$, $\rho_\ell := \inf\{t : \xi_t<\frac1\ell\}$, we may use standard results under Lipschitz assumptions). 
We fix $s=0$; the Markov property implies that is sufficient to prove the assertions of lemma 2.1 in this case. 

2)~ It is sufficient to consider case $\si(\cdot)\equiv 1$ only: the deterministic time change 
$$
A(t) \;:=\; \int_0^t \si^2(v)\, dv  \quad ,\quad  \tau(t) \;:=\; \inf\{ v>0 : A(v)>t \} 
$$ 
associates to $\,\wt\xi := (\xi_{\tau(t)})_{t\ge 0}$  a Brownian motion $\wt W$ such that $(*,s,u)$ is transformed to 
\beqq\label{decayeq2}
d\wt\xi_t \;=\; -\wt\beta(t)\, \wt\xi_t\, dt  \;+\; \sqrt{\wt\xi_t^+}\; d\wt W_t  \;,\quad 
t\ge 0 \;,\quad \wt\xi_0=u  
\eeqq  
where $\wt\beta(t)=\beta(\tau(t))$ is continuous, nonnegative and bounded; this follows from 
$$
\xi_{\tau(t)} - \xi_{\tau(s)} \;=\; - \int_{\tau(s)}^{\tau(t)} \beta(v)\, \xi_v\, dv 
\;+\; \int_{\tau(s)}^{\tau(t)} \si(v)\,\sqrt{\xi_v^+}\; dW_v \;. 
$$ 
By assumption (\ref{localization1}) on upper and lower bounds for $\si(\cdot)$, we have positive constants $\underline{c},\ov{c}$ such that 
$\underline{c}(t-s) \le \tau(t)-\tau(s) \le \ov{c}(t-s)$. The time change $t\to\tau(t)$ being deterministic, it is sufficient to prove the assertion of lemma 2.1 in the form (\ref{decayeq2}). Hence we assume $\si(\cdot)\equiv 1$ from now on.

3)~ Consider first the unique strong solution $\wt m$ for the problem without drift 
\beqq\label{decayeq4}
d\wt m_t \;=\; \sqrt{\wt m_t^+}\; d\wt W_t  \;,\quad t\ge 0 \;,\quad \wt m_0=u 
\eeqq
and write $\wt\rho$ for the time where $\wt m$ enters $0$. We have $\wt\rho<\infty$  $\,P$-almost surely by [IW 89, p.~237]. 
$\wt m$ is a local martingale before time $\wt\rho$, and a nonnegative supermartingale on $[0,\infty)$. 
Introduce stopping times  $\wt\rho^{(\ell)}<\wt\rho$ and stopped  processes $\wt m^{(\ell)}$  
$$
\wt\rho^{(\ell)} \;:=\; \inf\{ t>0: \wt m_t \notin [\frac1\ell,\ell]\} \;, \quad 
\wt m^{(\ell)} \;:=\; ( \wt m _{t\wedge\wt\rho^{(\ell)}})_{t\ge 0}
$$
for $\ell$ large enough. Then  $\wt m^{(\ell)}$ is a uniformly integrable martingale, and (\ref{decayeq4}) yields 
\beqq\label{decayeq5}
E\left( \langle \wt m^{(\ell)}-\wt m^{(\ell)}_0 \rangle _T\right) \;=\;  E\left( \langle \wt m^{(\ell)}-\wt m^{(\ell)}_0 \rangle _{T\wedge\wt\rho^{(\ell)}}\right) \;=\; 
E\left( \int_0^{T\wedge\wt\rho^{(\ell)}} \wt m^{(\ell)}_t\;  dt \right) \;\le\; u\cdot T
\eeqq
for all $T$. By $\,\wt\rho^{(\ell)}\uparrow\wt\rho\,$ as $\ell\to\infty$ and by monotone convergence   
$$
E\left( \sup_{0\le t\le T\wedge\wt\rho} (\wt m_t)^2 \right) 
\;=\; \sup_{\ell\in\bbn}\; E\left(  \sup_{0\le t\le T}  (\wt m^{(\ell)}_t)^2 \right) \;.  
$$
Since $\wt m$ is absorbed in $0$ at time $\wt\rho$, this shows in combination with (\ref{decayeq5})
\beqq\label{decayeq6}
E\left( \sup_{0\le t\le T} (\wt m_t)^2 \right)  
\;\le\;  \sup_{\ell\in\bbn}\; 2\left( u^2 + E\left( \sup_{0\le t\le T} (\wt m^{(\ell)}_t-\wt m^{(\ell)}_0)^2 \right) \right) 
\;\le\; C\left( u^2 + u\cdot T \right)   
\eeqq
for some constant $C$, using [IW 89, p. 110]. On the right hand side of (\ref{decayeq6}), the term which is linear in $u$ arises as a consequence of the square root form of the diffusion coefficient. 

4)~ For the process $\wt\xi$ solving equation (\ref{decayeq2}), define  
$\wt\tau := \inf\{t>0:\wt\xi_t=0\}$. $\;\beta(\cdot)$ being nonnegative, 
we may apply [KS 91, pp. 291+293] which does not need Lipschitz conditions,   
and compare pathwise the solution of (\ref{decayeq2}) to the solution of (\ref{decayeq4}): we obtain first    
$$
\wt\xi_t \;\le\; \wt m_t \quad\mbox{for all $t\in [[0,\wt\tau[[$, $P$-almost surely}     
$$ 
and second  
$$
\wt\xi_t \;\le\; \wt m_t \quad\mbox{for all $t\in [0,\infty)$}   
$$
since $\wt\xi_t$ is absorbed in $0$ at time $\wt\tau$ whereas $\wt m$ stays nonnegative. Thus step 3) implies  $\wt\tau<\infty$ $P$-almost surely, and we get in combination with (\ref{decayeq6})
$$
E\left( \sup_{0\le t\le T} (\wt\xi_t)^2 \right) \;\le\; E\left( \sup_{0\le t\le T} (\wt m_t)^2 \right)  
\;\le\; C\left( u^2 + u\cdot T \right) \;. 
$$
This gives the bound in the asserted form for equation (\ref{decayeq2}). The proof is finished. \halmos\\

Similiar to Protter [P 05, p.\ 250],  we will consider for  $p=1,2$ and for  $0<T<\infty$ norms 
\beqq\label{protternorm}
\left\| X \right\|_{*pT} \;:=\;  \left( E (\; \sup_{t\in [0,T]} |X_s|^p \;) \right)^{\frac1p}   
\eeqq
for c\`adl\`ag processes $X=(X_t)_{t\ge 0}$. Such norms with $p=2$ appear  in lemma 2.1.\\

{\bf 2.2 Lemma: }  For $s\ge 0$, $i=1,2$ and $u_1, u_2 >0$, consider strong solutions  
$$
\xi^{(s,u_i)}  \quad\mbox{of equations $(*,s,u_i)$ driven by $W^{(i)}$}   
$$
where $W^{(i)}$, $i=1,2$, are independent: then the sum 
$$
\xi \;:=\; \xi^{(s,u_1)}+\xi^{(s,u_2)}
$$
is a weak solution $\xi$  to equation $(*,s,u_1+u_2)$. \\
 
{\bf Proof: } This is the well-known 'branching property' of the CIR diffusion. Write $\xi^{(i)}:=\xi^{(s,u_i)}$, and prepare a third standard Brownian motion $W^{(0)}$ independent of $W^{(i)}$, $i=1,2$. Then 
$$
d\wh W_t \;:=\; 1_{\{ \xi^{(1)}_t+\xi^{(2)}_t>0 \}} \sum_{i=1,2} \sqrt{\frac{ \xi^{(i)}_t }{ \xi^{(1)}_t+\xi^{(2)}_t }}\, dW^{(i)}_t \;+\; 1_{\{ \xi^{(1)}_t+\xi^{(2)}_t = 0 \}} dW^{(0)}_t
$$
is again a standard Brownian motion. Equations $(*,s,u_i)$ driven by $W^{(i)}$, $i=1,2$, show that $\xi=\xi^{(1)}+\xi^{(2)}$ is a solution for equation $(*,s,u_1+u_2)$ driven by $\wh W$. \halmos\\

{\bf 2.3 Lemma: } For $0<a<b<\infty$ arbitrary, fixed, consider the equation     
$$
d\wh\xi_t \;=\; -\beta(t)\, \wh\xi_{t^-}\, dt  \;+\; \int_{(a,b]} y\; \mu(dt,dy) 
\;+\; \si(t)\,\sqrt{(\wh\xi_{t^-})^+}\; dW_t  \;,\quad t\ge 0  \;,\;  \wh\xi_0=0
\leqno(+,(a,b])
$$ 
driven by $W$ and $\mu$, both independent. 

a)~  Pathwise uniqueness holds for the SDE $(+,(a,b])$. 

b)~  Arrange the point masses of $\mu$ on $(0,\infty){\times}(a,b]$ in increasing order of time  
$$
\mu\left( \,\cdot\, \cap (0,\infty){\times}(a,b] \right) \;=:\; \sum_{i=1}^\infty \epsilon_{(T_i,Y_i)} 
\quad\mbox{where}\;\; 0<T_1 <T_2<\ldots \;\;\mbox{and}\;\; \lim_{i\to\infty}T_i=\infty \;, 
$$
prepare Brownian motions $W^{(i)}$, $i\ge 1$, independent and independent of $\mu$, and strong solutions 
$$
\mbox{ $\xi^{(T_i,Y_i)}$ for SDE $(+,T_i,Y_i)$ driven by $W^{(i)}$,  $i\ge 1$} \;. 
$$
Then 
$$
\wh\xi  \;=\; \wh\xi ^{(a,b]}  \;:=\;   \sum_{i=1}^\infty  \xi^{(T_i,Y_i)}\; 1_{[[T_i,\infty[[}   
$$
is a weak solution to $(+,(a,b])$. 

c)~ Equation $(+,(a,b])$ has a unique strong solution $\wh \xi = F(W,\mu)$. 

d)~ Using the norms defined in (\ref{protternorm}), the solution of $(+,(a,b])$ satisfies     
\beqq\label{bound1} 
\left\|  \wh\xi ^{(a,b]} \right\|_{*1T} 
\;\;\le\;\;  C\; T\sqrt{1+T}\; \int_{(a,b]}\sqrt{y+y^2}\; \nu(dy)    
\eeqq
\beqq\label{bound1-bis}
\left\|  \wh\xi ^{(a,b]} \right\|^2_{*2T} 
\;\;\le\;\;  C\; T (1+T)^2\;\; \nu( (a,b] ) \int_{(a,b]} (y+y^2)\; \nu(dy)    
\eeqq
for suitable constants $C$ which only depend on upper and lower bounds for $\si(\cdot)$ and on upper bounds for $\wt a(\cdot)$, according to (\ref{localization1}).\\

{\bf Proof: } 1)~ Extend the nonnegative processes $\xi^{(s,u)}$ of lemma 2.1 with starting time $s>0$ to the full time interval $[0,\infty)$ by $\xi_t^{(s,u)}\equiv 0$ for $0\le t<s$.  For any collection of points $(s_i,u_i)$ where $u_i>0$ and $0<s_1<s_2<\ldots$ such that $\lim\limits_{i\to\infty} s_i = \infty$, write $\xi^{(s_i,u_i)} = F(W^i,s_i,u_i)$ for strong solutions of $(+,s_i,u_i)$ driven by independent Brownian motions $W^i$. Fix $0<T<\infty$. As a consequence of lemma 2.1,  and by the norm properties of  (\ref{protternorm}), we have 
\beqq\label{bound2}
 \left\| \sum_{i=1}^\ell  \xi^{(s_i,u_i)}  \right\|_{*1T}  \le\;\;  \left\| \sum_{i=1}^\ell  \xi^{(s_i,u_i)}  \right\|_{*2T}
\le\;\; \sqrt{C(1+T)}\; \sum_{i=1}^\ell \sqrt{u_i+u_i^2} \;. 
\eeqq

2)~ With probability 1 we can arrange the points of $\mu(dt,dy)$ on $(0,\infty){\times}(a,b]$ as 
$$
(T_i,Y_i) \;,\; i\ge 1 \;,\quad\mbox{such that}\;\; 0<T_1 <T_2<\ldots \;\;\mbox{and}\;\; \lim_{i\to\infty}T_i=\infty \;.
$$
In this representation,  $\;(Y_i)_i$ are iid with values in $(a,b]$ and with law $\frac{\nu(\,\bullet\,\cap(a,b])}{\nu((a,b])}$. The sequences $(T_i)_i$ and $(Y_i)_i$ are independent. Similiar representations have been used in [HJ 93]. The number $N_T := \mu(\, (0,T]{\times}(a,b] \,)$ 
of points $(T_i,Y_i)$ with $T_i\le T$ is Poisson distributed with parameter 
$$
\nu((a,b])\; \int_0^T \wt a(s)\,ds \;\;\le\;\; \wt C\; T\; \nu((a,b])   
$$
by definition of Poisson random measure  with intensity  $\wt a(s) ds\, \nu(dy)$, and by assumption on $\wt a(\cdot)$. 

3)~ For PRM $\mu$ of step 2) independent of the BM's $W^i$, $i\ge 1$ of step 1), 
we may read inequality (\ref{bound2}) with $p=2$  and both sides squared as an inequality for conditional expectations given $\mu$: 
\beqq\label{bound3}
E \left(  \left[ \sup_{t\in [0,T]} \sum_{i=1}^{N_T} \xi_t^{(T_i,Y_i)} \right]^2 \mid \calf_\mu \right) 
\;\;\le\;\;   C(1+T)\;  \left[ \sum_{i=1}^{N_T} \sqrt{Y_i+Y_i^2}  \right]^2       
\eeqq 
where  
$$
\calf_\mu \;=\;  \si\left\{\, \mu\left( (0,t]{\otimes}B \right) : t\ge 0 , B\in\calb(\bbr) \,\right\}
$$
 denotes the $\si$-field generated by $\mu$ and where 
$$
\sum_{i=1}^{N_T} \sqrt{Y_i+Y_i^2}  \;=\;  \int_0^T\int_{(a,b]} \sqrt{y + y^2}\;\; \mu(ds,dy) \;. 
$$
Using the elementary inequality $(a+b)^2\le 2(a^2+b^2)$, we get 
\beao
\left[ \sum_{i=1}^{N_T} \sqrt{Y_i+Y_i^2}  \right]^2     
\;\le\;\;  2\left( \left[  \int_0^T\int_{(a,b]} \sqrt{y + y^2}\;\; \wt\mu(ds,dy) \right]^2  +  \left[  \int_0^T\int_{(a,b]} \sqrt{y + y^2}\; \wh\mu(ds,dy) \right]^2\right) \;.
 \eeao
Here the second term on the right hand side is deterministic since $\wh\mu(ds,dy)=\wt a(s)ds\,\nu(dy)$; using bounds for $\wt a(\cdot)$ according to (\ref{localization1}) and Jensen inequality, this term  is smaller than   
 $$
C\, T^2\left(\nu((a,b])\right)^2\; \left[ \int_{(a,b]} \sqrt{y + y^2}\, \frac{\nu(dy)}{\nu((a,b])} \right]^2  
  \;\;\le\;\;  C\, T^2\; \nu((a,b])\; \int_{(a,b]} (y + y^2)\, \nu(dy)  \;.
 $$
The first term on the right hand side is the square of a martingale  and thus has expectation   
$$
\int_0^T\int_{(a,b]} (y + y^2)\; \wh\mu(ds,dy)  \;\;\le\;\; C\,T\; \nu((a,b])\; \int_{(a,b]} (y + y^2)\, \nu(dy)   
$$
where we use again (\ref{localization1}). Putting all this together, we obtain 
\beqq\label{bound1-bis-provisorisch}
E \left(  \left[ \sup_{t\in [0,T]} \sum_{i=1}^{N_T} \xi_t^{(T_i,Y_i)} \right]^2  \right)  \;\;\le\;\;  C\, T(1+T)^2\; \nu((a,b])\; \int_{(a,b]} (y + y^2)\, \nu(dy)  
\eeqq
for some constant $C$ which only depends on upper and lower bounds for $\si(\cdot)$ (through (\ref{bound3}) and lemma 2.1) and upper bounds for $\wt a(\cdot)$. The right hand side of (\ref{bound1-bis-provisorisch}) corresponds to  (\ref{bound1-bis}). 
 
4)~To prove (\ref{bound1}), we read (\ref{bound2}) with $p=1$ conditionally on $\mu$ as
\beqq\label{bound3-bis}
E \left( \sup_{t\in [0,T]} \sum_{i=1}^{N_T} \xi_t^{(T_i,Y_i)}  \mid \calf_\mu \right) 
\;\;\le\;\;  \sqrt{ C(1+T) }\;  \sum_{i=1}^{N_T} \sqrt{Y_i+Y_i^2}   \;.   
\eeqq
Again, by the bounds on $\wt a(\cdot)$, we have  
$$
 E\left(  \int_0^T\int_{(a,b]} \sqrt{y + y^2}\; \mu(ds,dy)  \right) \;\;=\;\; \int_0^T\int_{(a,b]} \sqrt{y + y^2}\; \wh\mu(ds,dy)
  \;\;\le\;\; C\, T\; \int_{(a,b]} \sqrt{y + y^2}\; \nu(dy)
$$
and thus 
 \beqq\label{bound1-ter}
 E \left( \sup_{t\in [0,T]} \sum_{i=1}^{N_T} \xi_t^{(T_i,Y_i)}  \right)  \;\;\le\;\;  C\; T\sqrt{1+T}\; \int_{(a,b]}\sqrt{y+y^2}\; \nu(dy)    
\eeqq
where the right hand side of (\ref{bound1-ter}) corresponds to  (\ref{bound1}).

5)~ We show that with notations as above, the process  
$$ 
\wh\xi  \;:=\; \sum\limits_{i\ge 1} \xi^{(T_i,Y_i)} \;=\;   \sum\limits_{i\ge 1} \xi^{(T_i,Y_i)}1_{[[T_i,\infty[[}
$$ 
provides a weak solution to equation $(+,(a,b])$. To see this, prepare as in the proof of lemma 2.2 an additional standard Brownian motion $W^{(0)}$ independent of $W^{(i)}$, $i\in\bbn$, and consider 
$$
d\wh W_t \;:=\; 1_{\{ \wh\xi _{t-}>0 \}} \sum_{i=1}^\infty \sqrt{\frac{ \xi^{(T_iY_i)}_{t-} }{ \wh\xi _{t-} }}\, dW^{(i)}_t \;+\; 1_{\{ \wh\xi _{t-} = 0 \}} dW^{(0)}_t \;. 
$$
Then $\wh W$ is again a standard Brownian motion, and we have  
$$
d\wh\xi_t \;=\; -\beta(t)\, \wh\xi_{t^-}\, dt  \;+\; \int_{(a,b]} y\; \mu(dt,dy) 
\;+\; \si(t)\,\sqrt{(\wh\xi_{t^-})^+}\; d\wh W_t  \;,\quad t\ge 0  \;,\;  \wh\xi_0=0
$$
by definition of the processes $\xi^{(T_i,Y_i)}$, $i\in\bbn$. Note that in restriction to every interval $[[0,T_n[[$, the sums reduce to finite sums. 
 
 6)~ We finish the proof of lemma 2.1. First, pathwise uniqueness holds for solutions of equation $(+,(a,b])$ by Yamada-Watanabe, as in [H 08]. Second, we have constructed a weak solution to  $(+,(a,b])$ in the preceding step 3). Hence, as in Ikeda and Watanabe ([IW 89, theorem IV.1.1]), there is a unique strong solution for equation  $(+,(a,b])$ which has the form of a functional $F(W,\mu)$ of the driving pair $(W,\mu)$. We have to justify the last assertion:  
 
 Inspecting the proof  [IW 89, p.\ 163-166]  for  classical SDE that path\-wise uniqueness and existence of a weak solution together imply existence of a strong solution, functional of the driving Brownian path ([IW 89, theorem IV.1.1]),  we see that this proof carries over to SDE with jumps where the solution process is c\`adl\`ag. In fact, the proof is based only on the following: 
 i) c\`adl\`ag driving processes with stationary and independent increments; 
 ii)  transition probabilites between path spaces which are Polish. 
 So we can replace  $C$ used there (as path space for the solutions) by $D$, and $C_0=\{\al\in C:\al(0)=0\}$ used there (as path space for the driving process) by $C_0{\times}D_0$ since $D_0=\{\al\in D:\al(0)=0\}$ is a path space for $\left(\int_0^t\int_{(0,\infty)}y\, \mu(ds,dy)\right)_{t\ge 0}$ under assumption (\ref{permanentcondition}).\halmos\\

Now we turn to the problem 
$$
d\wh\xi_t \;=\; -\beta(t)\, \wh\xi_{t^-}\, dt  \;+\; \int_{(0,1]} y\; \mu(dt,dy) 
\;+\; \si(t)\,\sqrt{(\wh\xi_{t^-})^+}\; dW_t  \;,\quad t\ge 0  \;,\;  \wh\xi_0=0 \;.
\leqno(+,(0,1])
$$ 
For measures $\nu$ where $\nu((0,1])$ is finite, equation $(+,(0,1])$ has  finitely many jumps in finite time intervals, thus we can proceed as in lemma 2.3 to construct a weak solution and --~thanks to pathwise uniqueness~-- a unique strong solution $\wh\xi = F(W,\mu)$. It remains to consider measures $\nu$ having $\nu((0,1])=\infty$. If we think of $a$ tending to $0$ in lemma 2.3 and compare the structure of bounds (\ref{bound1}) and (\ref{bound1-bis}) for $a\downarrow 0$, the restrictive condition (\ref{restrictivecondition})
$$ 
\int_{(0,\infty)} (\sqrt{y}\wedge 1)\; \nu(dy) \;\;<\;\;  \infty   
$$ 
seems unavoidable  in this case. We start with a weak solution for $(+,(0,1])$ under (\ref{restrictivecondition}). \\

{\bf 2.4 Lemma: } Assume (\ref{restrictivecondition}) and $\nu((0,1])=\infty$. Fix some sequence $\delta_n \downarrow 0$ such that  
\beqq\label{choice}
0 \;\;<\;\; \int_{(0,\delta_n]} \sqrt{y}\; \nu(dy) \;\;<\;\;  4^{-n} \quad,\quad n\ge 1 \;, 
\eeqq 
and $\delta_0=1$. Prepare BM's $W^{(n)}$, $n\ge 1$, independent  and independent from $\mu$, together with 
$$
\mbox{$\wh\xi^{(\delta_n,\delta_{n-1}]}\;$  strong solution to $\,(+,(\delta_n,\delta_{n-1}])\,$  driven by $(W^{(n)},\mu)$}  
$$
for $n\ge 1$. Define for $N\ge 1$  
$$
\wh\xi^{(N)}  \;\;:=\;\;  \sum_{n=1}^N \wh\xi^{(\delta_n,\delta_{n-1}]}  \;. 
$$
Then the paths of $\wh\xi^{(N)}$ converge uniformly on compacts as $N\to\infty$ to the paths of  $\wh\xi=\wh\xi ^{(0,1]}$
\beqq\label{limitprocess}
\wh\xi^{(0,1]} \;\;:=\;\;  \sup_{N\ge 1} \wh\xi^{(N)} \;\;=\;\; \sum_{n=1}^\infty \wh\xi^{(\delta_n,\delta_{n-1}]}  
\eeqq
where the process defined by (\ref{limitprocess}) is c\`adl\`ag, provides a weak solution to the problem $\,(+,(0,1])\,$,  and satisfies the bound 
\beqq\label{bound5}
\left\|  \wh\xi ^{(0,1]} \right\|_{*1T} \;\;\le\;\;  C\; T\sqrt{1+T}\; \int_{(0,1]}\sqrt{y+y^2}\; \nu(dy)   
\eeqq
for some constant $C$ which depends only on bounds for the functions in (\ref{localization1}). \\

{\bf Proof: } We assume that the total mass $\nu((0,1])$ equals $+\infty$, and fix $0<T<\infty$. By choice (\ref{choice}) of $(\delta_n)_n$, by   
the bound (\ref{bound1}) in lemma 2.3 combined with Chebychev inequality, we have 
$$
P\left( \sup_{0\le t\le T}\; \wh\xi_t^{(\delta_n,\delta_{n-1}]}  \;>\;  2^{-n}  \right) \;\;\le\;\;  \wt C\, T\sqrt{1+T}\;\; 2^{-n}   
$$
for all $n\ge 1$. Using Borel Cantelli, we see that  $P$-almost surely  the paths of 
$$
\wh\xi^{(N)} \;\;=\;\; \sum_{n=1}^N \wh\xi^{(\delta_n,\delta_{n-1}]}  
$$
converge uniformly on $[0,T]$ as $N\to\infty$ to a finite limit 
$$
\wh\xi \;\;:=\;\;  \sup_{N\ge 1} \wh\xi^{(N)} \;. 
$$
A uniform limit of c\`adl\`ag processes is again c\`adl\`ag. By monotone convergence, and applying the bounds (\ref{bound1}) of lemma 2.3 to all processes $\wh\xi^{(N)}$, $N\ge 1$,  we obtain for the limit  (\ref{limitprocess}) the bound  
$$
E\left( \,\sup\limits_{0\le t\le T}\, \wh\xi_t\,  \right)  
\;\;\le\;\;  C\; T\sqrt{1+T}\; \int_{(0,1]}\sqrt{y+y^2}\; \nu(dy)   
$$
thanks to assumption (\ref{restrictivecondition}). This is the bound (\ref{bound5}). 

2)~ We prove that the limit  in (\ref{limitprocess})  is indeed a weak solution to equation  $(+,(0,1])$.  Due to the bound (\ref{bound5}), the process 
$$
\sum_{n=1}^\infty   \int_0^t \sqrt{\wh\xi^{(\delta_n,\delta_{n-1}]}_{s-}}\; dW^{(n)}_s \;,\quad t\ge 0
$$
is a well defined martingale whose angle brackett 
$$
\sum_{n=1}^\infty   \int_0^t  \wh\xi^{(\delta_n,\delta_{n-1}]}_{s-}\; ds  \;=\;  \int_0^t  \wh\xi_{s-}\; ds  \;,\quad t\ge 0  
$$
is integrable at every fixed time $t<\infty$. Prepare as in step 4) of the proof of lemma 2.3 another standard Brownian motion $W^{(0)}$ independent of $W^{(n)}$, $n\in\bbn$, and of $\mu$, and consider 
$$
d\wh W_t \;:=\; 1_{\{ \wh\xi _{t-}>0 \}} \sum_{n=1}^\infty \sqrt{\frac{ \wh\xi^{(\delta_n,\delta_{n-1}]}_{t-} }{ \wh\xi _{t-} }}\, dW^{(n)}_t \;+\; 1_{\{ \wh\xi _{t-} = 0 \}} dW^{(0)}_t \;. 
$$
Then $\wh W$ is again a standard Brownian motion, and we do have  
$$
d\wh\xi_t \;=\; -\beta(t)\, \wh\xi_{t^-}\, dt  \;+\; \int_{(0,1]} y\; \mu(dt,dy) 
\;+\; \si(t)\,\sqrt{(\wh\xi_{t^-})^+}\; d\wh W_t  \;,\quad t\ge 0  \;,\;  \wh\xi_0=0
$$
by definition of the processes $\wh\xi^{(\delta_n,\delta_{n-1}]}$, $n\in\bbn$:  
here we exploit convergence  $\wh\xi^{(N)}\uparrow\wh\xi$ a.s.\ uniformly on compacts as $N\to\infty$, proved in step 1), to show convergence of drift terms  
$$
\int_0^\cdot \beta(t)\, \wh\xi^{(N)}_{t^-}\, dt   \;\;\lra\;\;  \int_0^\cdot \beta(t)\, \wh\xi_{t^-}\, dt \quad,\quad N\to\infty
$$
and angle bracketts of martingale terms 
$$
\left\langle   \int_0^\cdot \si(t)\,\sqrt{(\wh\xi_{t^-})^+}\; d\wh W_t 
\;-\;  \sum_{n=1}^N \int_0^\cdot \si(t)\, \sqrt{ \wh\xi^{(\delta_n,\delta_{n-1}]}_t }\;  dW^{(n)}_t \right\rangle \;\;\lra\;\;  0 \quad,\quad N\to\infty 
$$
a.s.\ uniformly on compact time intervals. Since  $\mu$ induces positive jumps only which are summable, it is obvious that 
$$
\int_0^\cdot \int_{(\delta_N,1]}y\, \mu(ds,dy)   \;\;\lra\;\; \int_0^\cdot \int_{(0,1]} y\, \mu(ds,dy)  \quad,\quad N\to\infty
$$
a.s.\ uniformly on compact time intervals. The three last assertions together prove that the limit process $\wh\xi$ in (\ref{limitprocess}) provides  indeed a solution to the problem $\,(+,(0,1])\,$.  \halmos\\

{\bf 2.5 Lemma: } Under assumption (\ref{restrictivecondition}), for measures $\nu$ with $\nu((0,1])=\infty$, 

a)~ pathwise unicity holds for equation $(+,(0,1])$; 

b)~ equation $(+,(0,1])$ has a unique strong solution $\wh\xi = F(W,\mu)$. \\ 

{\bf Proof: } Exactly as in step 6) in the proof of lemma 2.3, based on existence of a weak solution which has been established in lemma 2.4.\halmos\\

{\bf 2.6 Proof of theorem 1.1: } 1)~ In order to solve equation $(\diamond)$ with starting point  $x_0\ge 0$, introduce BM's $W^{(i)}$, $i=1,2,3$, independent and independent  of  $\mu$, together with the following processes $\xi^{(i)}=F(W^{(i)},\mu)$: 

i)~ $\xi^{(1)}=F(W^{(1)},\mu)$ is a strong solution to  equation $(+,(0,1])$ according to lemma 2.5, driven by $W:=W^{(1)}$ and the small jumps of $\mu$; 

ii)~ $\xi^{(2)}=F(W^{(2)},\mu)$ is a strong solution to equation  
\beqq\label{bigjumpequation}
d \xi_t \;=\; -\beta(t)\,\xi_{t^-}\, dt  \;+\; \int_{(1,\infty)} y\; \mu(dt,dy) 
\;+\; \si(t)\,\sqrt{(\xi_{t^-})^+}\; dW_t  \;,\quad t\ge 0  \;,\;  \wh\xi_0=0  
\eeqq
driven by $W:=W^{(2)}$ and the big jumps of $\mu$ (since by assumption there are at most finitely many big jumps over finite time intervals, this strong solution is obtained as in steps 2)+5)+6) of lemma 2.3; an analogue of (\ref{bound3}) holds, whereas no finite moments and hence no moment bounds corresponding to (\ref{bound1}) or (\ref{bound1-bis}) are available for $\xi^{(2)}=F(W^{(2)},\mu)$ under our assumptions on $\nu$); 

iii)~ $\xi^{(3)}=F(W^{(3)},x_0)$ is a strong solution to the 'continuous' problem
\beqq\label{continuousequation}
d \xi_t \;=\; \left[ a(t)-\beta(t)\, \xi_t\right] dt  
\;+\; \si(t)\,\sqrt{( \xi_t)^+}\; d W_t  \;,\quad t\ge 0  \;,\;   \xi_0=x_0\ge 0  
\eeqq
with $x_0\ge 0$ the starting point of equation $(\diamond)$ 
(in case $x_0=0$, the process (\ref{continuousequation}) is identically $0$; in case $x_0>0$, a unique strong solution to (\ref{continuousequation})  exists as in lemma 2.1). 

2)~ From the processes defined in step 1), we construct a weak solution $(\wh \xi,\wh W)$ for equation $(\diamond)$ with starting point $x_0$ by
$$
\wh \xi \;:=\;  \sum_{i=1}^3 \xi^{(i)}  \quad,\quad 
d\wh W_t \;=\;  1_{\{\wh \xi_{t-}>0\}} \sum_{i=1}^3 \sqrt{\frac{\xi^{(i)} _{t-}}{\wh \xi_{t-}}} d W^{(i)} _t  
\;+\; 1_{\{\wh \xi_{t-}=0\}} dW^{(0)}_t \;. 
$$
Since pathwise uniqueness holds for equation $(\diamond)$ (extension of Yamada-Watanabe as in [H~08]) we conclude as in step 6) of the proof of lemma 2.3 (extension of [IW 89, theorem IV.1.1] to SDE with jumps) that there is a unique strong solution $\xi=F(W,\mu,x_0)$ for equation $(\diamond)$, functional of the driving pair $(W,\mu)$ and of the starting point $x_0$. Theorem 1.1 is now proved. \halmos\\

\section*{3. Proof of theorem 1.2} 

We give the proof of theorem 1.2  through a series of lemmata. Again, a localization argument shows that it is sufficient to consider the case where the assumptions of 1.2 are strengthened to   
\beqq\label{localization2}
\begin{array}{l}
\mbox{$\beta(\cdot)$ and $\si(\cdot)$ are  bounded away from both $0$ and $\infty$ on $[0,\infty)$}\\
\mbox{$\al(\cdot)$, $\wt\al(\cdot)$ are bounded on $[0,\infty)$} \;.  
\end{array}
\eeqq
We shall assume (\ref{localization2}) throughout this proof section. We start with asymptotics of  (\ref{3})+(\ref{3bis}) under condition  (\ref{localization2}), and use the notation  $a_t\asymp b_t$ in the sense of  $0<\liminf\limits_{t\to\infty} \frac{a_t}{b_t}\le \limsup\limits_{t\to\infty} \frac{a_t}{b_t}<\infty$.\\

{\bf 3.1 Lemma: } Under  (\ref{localization2}), for  $0\le s<t<\infty$ and $\gamma(s,t)$, $p(s,t)$, $B(s,t)$ defined in (\ref{3})+(\ref{3bis}):  
\beam
\gamma(s,t) &\sim& \frac{2}{\si^2(s)}\, \frac{1}{t-s}  \quad\quad\mbox{as}\;\; t\downarrow s \;,
\label{L1aii}\\  
p(s,t) &\sim&   \frac{2}{\si^2(s)}\, \frac{1}{t-s}  \quad\quad\mbox{as}\;\; t\downarrow s \;,  
\label{L1aiii}\\
p(s,t) 
&\asymp& 1 \quad\quad\mbox{as}\;\; t\uparrow \infty \;,
\label{L1bi}\\
\gamma(s,t) &\asymp& B(s,t)  \quad\mbox{as}\;\; t\uparrow \infty \;. 
\label{L1bii} 
\eeam \\[-7mm]
Under (\ref{localization2}),  for every choice of $0\le s<\infty$,  $B(s,t)$ decreases exponentially fast as $t\to\infty$. \\

{\bf Proof: } For $0\le s<t$ and $t$ decreasing to $s$, we have for $C(s,t)$ and $B(s,t)$ defined in (\ref{3})
$$
\left[ \frac{1}{t-s}C(s,t)\right] B(0,s)  \;\lra\;  \left[ \frac{\si^2(s)}{2}\, \exp\left( \int_0^s \beta(v)dv \right) \right] B(0,s) \;=\; \frac{\si^2(s)}{2} \quad\quad\mbox{as}\;\; t\downarrow s \;,
$$
which gives first (\ref{L1aii}) by definition of  $\gamma(s,t)$  in (\ref{3bis}), and second (\ref{L1aiii}) for $p(s,t)=\frac{\gamma(s,t)}{B(s,t)}$. 

For $s$ fixed and $t$ increasing to $\infty$, assumption (\ref{localization2}) guarantees 
$$
B(0,t)C(s,t) \;=\; \int_s^t  \frac{\si^2(v)}{2}\, \exp\left( -\int_v^t \beta(r)dr \right)\, dv  \;\;\asymp\;\; 1  \quad\quad\mbox{as}\;\; t\to\infty  
$$
which gives first (\ref{L1bi}) by definition of $p(s,t)$ in (\ref{3bis}), and second (\ref{L1bii}) for $\gamma(s,t)=B(s,t)p(s,t)$.  \halmos \\

{\bf 3.2 Lemma: } Consider $0\le s<t<\infty$. For $\la \to \Psi_{s,t}(\la)$ defined in (\ref{3ter}) we have 
\beqq\label{L2i}
\left(\frac{\partial}{\partial\la}\Psi_{s,t}\right)(0^+) \;=\; B(s,t) \;,\quad 
\left(\frac{\partial^2}{\partial\la^2}\Psi_{s,t}\right)(0^+) \;=\; - \frac{2B(s,t)}{p(s,t)}  
\eeqq
and for all $\la\ge 0$ 
\beqq\label{L2ii}
\lim_{t\downarrow s}\,\Psi_{s,t}(\la) \;=\;  \la \,,\quad \lim_{t\uparrow \infty}\,\Psi_{s,t}(\la) \;=\; 0 \;. 
\eeqq
Moreover, the family $\left(\Psi_{s,t}\right)_{0\le s<t<\infty}$ has a semigroup property under functional iteration: 
\beqq\label{L2iii}
\Psi_{t_1,t_2} \circ \Psi_{t_2,t_3} \;=\; \Psi_{t_1,t_3}
\quad\mbox{on $[0,\infty)$, for $0\le t_1 < t_2 < t_3 <\infty$}\;. 
\eeqq

\vskip0.5cm
{\bf Proof: } Calculating derivatives with respect to $\la$ in 
$$
\Psi_{s,t}(\la) \;=\; \gamma(s,t)\, \frac{\la}{p(s,t)+\la} \;,\quad\la\ge 0
$$ 
at $\la=0^+$ gives (\ref{L2i}). The limits (\ref{L2ii}) are an immediate consequence of lemma 3.1. 
We turn to the iteration property. Fix $0\le t_1 < t_2 < t_3 <\infty$. 
First, 
\beqq\label{L2B1}
p(t_1,t_2) + \gamma(t_2,t_3) \;=\; \frac{C(t_1,t_3)}{B(0,t_2)C(t_1,t_2)C(t_2,t_3)} 
\eeqq
follows from definition of $p(t_1,t_2)$ and $\gamma(t_2,t_3)$ since 
$$
\frac{1}{B(0,t_2)C(t_1,t_2)} + \frac{1}{B(0,t_2)C(t_2,t_3)} \;=\; 
\frac{1}{B(0,t_2)}\, \frac{C(t_1,t_2)+C(t_2,t_3)}{C(t_1,t_2)C(t_2,t_3)}  
$$
where the numerator is by definition additive. Using (\ref{L2B1}), the definition of $p(.,.)$ and $\gamma(.,.)$ 
gives 
\beqq\label{L2B2}
\frac{p(t_1,t_2)p(t_2,t_3)}{p(t_1,t_2) + \gamma(t_2,t_3)} \;=\;  p(t_1,t_3) \;,\quad 
\frac{\gamma(t_1,t_2)\gamma(t_2,t_3)}{p(t_1,t_2) + \gamma(t_2,t_3)} \;=\;  \gamma(t_1,t_3) \;. 
\eeqq
Using (\ref{L2B1})+(\ref{L2B2}), we calculate 
\beao
\Psi_{t_1,t_2} \circ \Psi_{t_2,t_3}(\la)   &=&    \Psi_{t_1,t_2} \left(\Psi_{t_2,t_3}(\la)\right)  
\;\;=\;\; \gamma(t_1,t_2)\; \frac{ \gamma(t_2,t_3)\frac{\la}{p(t_2,t_3)+\la} }{ p(t_1,t_2) + \gamma(t_2,t_3)\frac{\la}{p(t_2,t_3)+\la} } \\
&=&   \frac{ \gamma(t_1,t_2)\gamma(t_2,t_3)\;\la }{ p(t_1,t_2)p(t_2,t_3) + \la(p(t_1,t_2)+\gamma(t_2,t_3)) } \\[2mm]
&=&   \gamma(t_1,t_3)\frac{\la}{p(t_1,t_3)+\la}  \quad=\quad \Psi_{t_1,t_3}(\la)   
\eeao
which is the iteration property (\ref{L2iii}). In this last part of the proof, we follow closely the representation given in Hammer [H 06] for the classical CIR diffusion. Lemma 3.2 is proved. \halmos \\

The next step is to identify a semigroup  $\left(H_{s,t}(\cdot,\cdot)\right)_{0\le s<t<\infty}$ of transition probabilities on the state space  $([0,\infty),\calb([0,\infty))$ of $(\diamond)$ with the property 
\beqq\label{semigroupH}
H_{s,t}(u,B) \;=\; P\left( \xi ^{(s,u)}_t \in B \right) \quad,\quad 0\le s<t<\infty \,,\, u\ge 0 \,,\, B\in\calb([0,\infty))  
\eeqq
where $\,\xi ^{(s,u)}$ denotes the solution to equation $(*,s,u)$ in lemma 2.1. For the classical CIR diffusion, Poisson mixture representations of the transition probabilities have been considered in Chaleyat-Maurel and Genon-Catalot [CG 06, (78)+(82)]; see also G\"oing-Jaeschke and Yor [GY 03] for some context.\\

{\bf 3.3 Lemma : } a)~  For $0\le s<t<\infty$, introduce a transition  probability $H_{s,t}(\cdot,\cdot)$ 
$$
H_{s,t}(y,\cdot)  \;:=\;  e^{\, -y\,\gamma(s,t)}\,\epsilon_0(\cdot) +  
\sum_{m=1}^\infty  \frac{ e^{\, -y\,\gamma(s,t)}\,(\,y\,\gamma(s,t))^m }{m\,!}\, \Gamma(m,p(s,t))   \quad,\quad y\ge 0 
$$
on $([0,\infty),\calb([0,\infty))$. For $y>0$, $\,H_{s,t}(y,\cdot)$ is a Poisson mixture of Gamma laws; for $y=0$,  
$H_{s,t}(0,\cdot)$ reduces to Dirac measure $\epsilon_0$ sitting in $0$. For all $y\ge 0$, the Laplace transform 
$$
\vp_{s,t}(y;\la) \;:=\; \int_0^\infty e^{-\la x}\, H_{s,t}(y,dx) \quad,\; \la\ge 0
$$
of $H_{s,t}(y,\cdot)$ is given by Barra [B 71, Ch. VII.1] as  
$$
\vp_{s,t}(y;\la) 
\;=\;  \exp\left( -y\;\gamma(s,t)\;\frac{\la}{p(s,t)+\la} \right)  \;=\;  \exp\left( -y\,\Psi_{s,t}(\la) \right) \;.
$$
 
b)~ For fixed $s$ and  $y>0$, under (\ref{localization2}), mean and variance of $H_{s,t}(y,\cdot)$   
\beqq\label{L3i}
\int_0^\infty x\, H_{s,t}(y,dx) \;=\; y\,B(s,t) \;,\quad 
\int_0^\infty (x-yB(s,t))^2\, H_{s,t}(y,dx) \;=\; y\,\frac{2\,B(s,t)}{p(s,t)}  
\eeqq
vanish exponentially fast when $t$ tends to $\infty$. 

\vskip0.2cm
c)~ $(H_{s,t}(\cdot,\cdot))_{0\le s<t<\infty}$ is a semigroup of transition probabilities on $\left( [0,\infty) , \calb([0,\infty)) \right)$. 

\vskip0.2cm
d)~ Write $\calc^\infty_\kappa$ for the class of $\calc^\infty$-functions having compact support contained in $(0,\infty)$. 
In restriction to $\left( (0,\infty) , \calb((0,\infty)) \right)$, $\;(H_{s,t}(\cdot,\cdot))_{0\le s<t<\infty}$ acts as a 
submarkovian semigroup with the property  
$$
\lim_{t\downarrow s}\frac{H_{s,t}f-f}{t-s}(y) \;=\; (\cala_s f)(y) 
\;:=\; - \beta(s)\, y\,  f'(y) + \frac{\si^2(s)}{2}\, y\, f''(y) \;, \quad f\in\calc^\infty_\kappa \;,\; y>0  
$$
for every $s$. 

e)~ As a consequence of d), the semigroup $(H_{s,t}(\cdot,\cdot))_{s<t}$ on $\left( [0,\infty) , \calb([0,\infty)) \right)$ corresponds to the system of SDE  $\left\{ (*,s,y) : 0\le s<\infty , y\ge 0 \right\}$ with absorption at $0$ considered in lemma 2.1. \\

{\bf Proof: } 1) Assertion a) follows from [B 71, p.\ 82]; e) is immediate from d); b) comes from derivatives of 
$\vp_{s,t}(y;\cdot)$ at $0^+$, using lemma 3.1 and (\ref{L2i}) in lemma 3.2. We have to prove c)+d). 

2) To $\la> 0$ we associate $z = e^{-\la}$, $0<z\le 1$,  and define for $0\le s<t<\infty$ and $y\ge 0$ 
$$
g_{s,t}(y;z)  \;:=\;  \vp_{s,t}(y;-\log(z))  \;=\; \int_{[0,\infty)} z^x\, H_{s,t}(y,dx)   \;,\quad 0<z\le 1  \;.
$$ 
Then the function $\;g_{s,t}(y;\cdot)$ is strictly increasing on  $(0,1]$, takes the value $g_{s,t}(y;1)=1$ at $z=1$, and satisfies by part a) of the lemma  
$$
g_{s,t}(y;\cdot) \;=\; \left( g_{s,t}(1;\cdot) \right)^y  \quad\mbox{for $0\le s<t<\infty$ and $y\ge 0$} \;. 
$$
It is easy to see that the set of functions 
$$
g_{s,t}(1;\cdot) \;:\;  (0,1]\to (0,1]   \quad,\quad   0\le s<t<\infty
$$
has a functional iteration property 
\beqq\label{L3iii}
g_{t_1,t_2}(1;\cdot) \circ g_{t_2,t_3}(1;\cdot) \;=\; g_{t_1,t_3}(1;\cdot)
\quad\mbox{for $0\le t_1 < t_2 < t_3 <\infty$}  
\eeqq
similiar to the iterates of the generating function in classical branching process theory (cf.\ Athreya and Ney [AN 72]): this is a consequence of (\ref{L2iii}) in lemma 3.2 since by definition  
$$
g_{t_1,t_2}(1;g_{t_2,t_3}(1;z)) \;=\; \vp_{t_1,t_2}(1;-\log \vp_{t_2,t_3}(1;\la)) \;=\;
e^{ -\Psi_{t_1,t_2}(\Psi_{t_2,t_3}(\la)) } \;. 
$$
From (\ref{L2ii}) in lemma 3.2 we have  immediately  
\beqq\label{L3ii}
\lim_{t\downarrow s}\, g_{s,t}(1;z) \;=\; z \quad,\quad \lim_{t\uparrow\infty}\, g_{s,t}(1;z) \;=\; 1 
\quad,\quad 0<z\le 1 \;. 
\eeqq

3)~ We prove the semigroup property of $\left(H_{s,t}(\cdot,\cdot)\right)_{0\le s<t<\infty}$ on $([0,\infty),\calb([0,\infty))$, and thus c): consider arbitrary points $0\le t_1<t_2<t_3<\infty$, then for $z\in (0,1]$
\beao
&&\int_{[0,\infty)}\int_{[0,\infty)} H_{t_1,t_2}(y,du)\, H_{t_2,t_3}(u,dx)\, z^x \;\;=\;\; 
\int_{[0,\infty)} H_{t_1,t_2}(y,du)\, g_{t_2,t_3}(u;z) \\
&&=\quad \int_{[0,\infty)} H_{t_1,t_2}(y,du)\, [g_{t_2,t_3}(1;z)]^u  \;\;=\;\;  g_{t_1,t_2}(y; g_{t_2,t_3}(1;z)) \\
&&=\quad \left[g_{t_1,t_2}(1; g_{t_2,t_3}(1;z) )\right]^y \;\;=\;\; \left[g_{t_1,t_3}(1;z)\right]^y 
\;\;=\;\; g_{t_1,t_3}(y;z)  \;\;=\;\; \int_{[0,\infty)} H_{t_1,t_3}(y,dx)\, z^x  \;. 
\eeao

4)~ To determine the generator of $(H_{s,t}(\cdot,\cdot))_{s<t}$ at $t=s$ when the semigroup is restricted  $(0,\infty)$, we give a heuristic argument. Consider $y>0$. 
By $(\ref{L3i})$ combined with $(\ref{L1aiii})$ and by definition of $B(s,t)$ in (\ref{3}), the laws $H_{s,t}(y,\cdot)$ are concentrated as $t\downarrow s$ on small neighbourhoods 
around their mean value $y\,B(s,t)$, which tends to $y$. Thus we have for $f\in\calc^\infty_\kappa$ as $t\downarrow s$  
\beao
\left( H_{s,t}f-f \right)(y)  &=& \int H_{s,t}(y,dx)\, [f(x)-f(y)] \\ 
&=& [f(yB(s,t))-f(y)] \;\;+\;\;    f'(yB(s,t)) \int H_{s,t}(y,dx)\, [x-yB(s,t)] \\
&&+\quad   \frac12f''(yB(s,t)) \int H_{s,t}(y,dx)\, [x-yB(s,t)]^2 \;\;+\;\; \ldots \;. 
\eeao
On the right hand side of the last equation, the second term is $0$ by $(\ref{L3i})$, the third  one equals   
$$
\frac{y\,B(s,t)}{p(s,t)}\,f''(yB(s,t)) \;\;\sim\;\; (t-s)\,\frac{\si^2(s)}{2}\, y\, f''(y)  
\quad\quad\mbox{as}\quad t\downarrow s 
$$
by (\ref{L3i})+(\ref{L1aiii}), and the first term behaves as  
$$
f'(y)\,y\, [B(s,t)-1]  \;\;\sim\;\;  -(t-s)\,\beta(s)\,y\, f'(y)  
\quad\quad\mbox{as}\quad t\downarrow s \;. 
$$
Together, this shows for $0\le s<t<\infty$ and for $y>0$
$$
\lim_{t\downarrow s}\,\frac{H_{s,t}f-f}{t-s}(y) \;=\; (\cala_s f)(y) 
$$
with $\cala_s$ as defined in d). This gives d), and concludes the proof of lemma 3.3. \halmos\\

{\bf 3.4 Lemma: } For $0\le s<t<\infty$, we have 
$$
\exp\left( - \int_s^t \frac{\si^2(v)}{2}\, \Psi_{v,t}(\la)\, dv\, \right) \;=\; \left( 1 + \frac{\la}{p(s,t)}\right) ^{-1} 
\;,\quad \la\ge 0 \;. 
$$
This is the Laplace transform of an exponential law with parameter $p(s,t)$.\\

{\bf Proof: } Taking first the log of both sides and deriving then with respect to $s$ for fixed $t$, 
we obtain on the left hand side and on the right hand side the same derivative  
$ 
\frac{\si^2(s)}{2}\, \Psi_{s,t}(\la) \;;  
$ 
when calculating the derivative on the right hand side, we make use of 
$$
\frac{\partial}{\partial s}\;\frac{1}{p(s,t)} \;=\; B(0,t)\, \frac{\partial}{\partial s}C(s,t) 
\;=\; -B(s,t)\,\frac{\si^2(s)}{2} \;. \eqno\Box
$$

\vskip0.8cm
{\bf 3.5 Proof of remark 1.3: } This is an immediate consequence of  lemma 3.4.\halmos\\

The next step is to specify the semigroup  $\left( J_{s,t}(\cdot,\cdot)\right)_{0\le s<t<\infty}$ of a CIR diffusion (\ref{continuousequation})
$$
d \xi_t \;=\; \left[ a(t)-\beta(t)\, \xi_t\right] dt  
\;+\; \si(t)\,\sqrt{( \xi_t)^+}\; d W_t  \;,\quad t\ge 0  
$$
with time dependent coefficients. We start with \\

{\bf 3.6 Lemma: } For given $\al(\cdot)$ nonnegative, right continuous and bounded, there is a family of laws $\left(I_{s,t}\right)_{0\le s<t<\infty}$ on $\left( [0,\infty),\calb([0,\infty))\right)$ with the 'skew convolution'  (Li [L 93]) property 
\beqq\label{L6ii}
I_{t_1,t_2}H_{t_2,t_3} * I_{t_2,t_3} \;=\; I_{t_1,t_3}  \quad\mbox{for all $0\le t_1<t_2<t_3<\infty$}  
\eeqq
relative to the semigroup  $\left(H_{s,t}\right)_{0\le s<t<\infty}$ of  lemma 3.3, and with Laplace transforms 
\beqq\label{L6i}
\quad\int_0^\infty e^{-\la x}\, I_{s,t}(dx)  \;\;=\;\;  
\exp\left( -\int_s^t \frac{\si^2(v)}{2} \al(v)\, \Psi_{v,t}(\la)\, dv \right) \;\;,\quad 0\le s<t<\infty \;,\; \la\ge 0 \;. 
\eeqq

\vskip0.8cm
{\bf Proof: } 1)~ Note first  that for constant $\al>0$ and $t_1<t_2<t_3$,   the law 
$$
\Gamma(\al,p(t_1,t_2))H_{t_2,t_3} \;=\; \int_{[0,\infty)} \Gamma(\al,p(t_1,t_2))(dy)\, H_{t_2,t_3}(y,\cdot)  
$$
has by 3.3 and 1.3 and (\ref{L2iii}) the Laplace transform 
$$
\la \;\lra\; \exp\left( -\int_{t_1}^{t_2} \frac{\si^2(v)}{2} \al\, \Psi_{v,t_2}(\Psi_{t_2,t_3}(\la))\, dv \right)  
\;=\; \exp\left( -\int_{t_1}^{t_2} \frac{\si^2(v)}{2} \al\, \Psi_{v,t_3}(\la)\, dv \right) \;. 
$$

2)~ Fix $0\le s<t<\infty$. For  the grid of points $t_{n,j}:=j2^{-n}$, $j\in\bbn_0$, 
introduce step functions  
$$
\al_n(v) \;:=\; \sum_{j=1}^\infty 1_{[t_{n,j-1},t_{n,j})}(v)\; \al(t_{n,j}) \;,\quad v\ge 0
$$
which -- by right continuity of $\al(\cdot)$ -- converge to $\al(\cdot)$ as $\nto$. For  $j\ge 1$ such that 
intervals $[t_{n,j-1},t_{n,j})$ and $(s,t]$ intersect, write $Q_{n,j}$ for the 
law on $[0,\infty)$ having Laplace transform 
$$
\la  \;\lra\;  \exp\left( - \int_{s\vee t_{n,j-1}}^{t\wedge t_{n,j}} \frac{\si^2(v)}{2}\, \al_n(v)\, \Psi_{v,t}(\la)\, dv\,\right) 
$$
as in step 1). Convolution of the $Q_{n,j}$ 
defines a law $Q_n$ on $[0,\infty)$ whose Laplace transform is 
\beqq\label{L6Bii}
\la \;\lra\;  \exp\left( - \int_s^t \frac{\si^2(v)}{2}\, \al_n(v)\, \Psi_{v,t}(\la)\,dv \, \right)  \;. 
\eeqq
Now the right hand side of (\ref{L6i}) viewed as function of $\la$ approaches the limit $1$ as $\la\downarrow 0$, and is the pointwise limit of Laplace transforms (\ref{L6Bii}) as $\nto$. Hence the function (\ref{L6i})  itself (Feller [F 71, p.\ 431]) is Laplace transform of some probability measure $I_{s,t}$ on $[0,\infty)$. 
 
3)~ Rephrase steps 1)+2) as follows: $I_{s,t}$ arises as weak limit of convolutions 
$$
\begin{array}{l}
\quad\Gamma\left(\, \al(r_{n,1}) \,,\, p(r_{n,0},r_{n,1}) \,\right) H_{r_{n,1} , t}   \;*\;\ldots\;*\; 
\Gamma\left(\, \al(r_{n,j}) \,,\, p(r_{n,j-1},r_{n,j}) \,\right) H_{r_{n,j} , t}  \;*\;\ldots  \\
\quad*\;  \Gamma\left(\, \al(r_{n,\ell(n)-1}) \,,\, p(r_{n,\ell(n)-2},r_{n,\ell(n)-1}) \,\right)  H_{r_{n,\ell(n)-1} , t}  
\;*\;  \Gamma\left(\, \al(r_{n,\ell(n)}) \,,\, p(r_{n,\ell(n)-1},r_{n,\ell(n)}) \,\right) 
\end{array}\\[2mm]
$$
along grids $s=r_{n,0}<r_{n,1}<\ldots<r_{n,l(n)-1}<r_{n,l(n)}=t$
in $(s,t]$ with mesh tending to $0$, as $\nto$, or equivalently as weak limit  of laws 
\beqq\label{convolutioninterpretation1}
\call\left(   \sum_{j=1}^{\ell(n)}  \xi^{(r_{n,j},U_{n,j})}_t  \right)  \quad,\quad \nto  
\eeqq
where at every stage $n$ of the asymptotics,  $U_{n,j} \,\sim\, \Gamma\left(\, \al(r_{n,j}) \,,\, p(r_{n,j-1},r_{n,j}) \,\right)$ are independent random variables and  $\xi^{(r_{n,j},U_{n,j})}$ independent processes solving $(*,r_{n,j},U_{n,j})$ of lemma 2.1, for $1\le j\le \ell(n)$. Note also that with high probability, the random variable $U_{n,j}$ takes its values in small neighbourhoods of $\frac{\si^2(r_{n,j})}{2}\,  \al(r_{n,j})\, (r_{n,j}-r_{n,j-1}) $, as a consequence of (\ref{L1aiii}) in lemma 3.1. 
From interpretation  (\ref{convolutioninterpretation1}) of $I_{s,t}$, the skew convolution property  (\ref{L6ii}) is obvious.   \halmos\\

The interpretation of $I_{s,t}$ in terms of limits of laws (\ref{convolutioninterpretation1}) reappears in the following lemma: \\

{\bf 3.7 Lemma: } For given $\al(\cdot)$ nonnegative,  right continuous and  bounded, the semigroup of transition probabilities associated to a CIR diffusion with time dependent coefficients 
\beqq\label{continuousequationspecialrepresentation}
d \xi_t \;=\; \left[ \frac{\si^2(t)}{2}\al(t)-\beta(t)\, \xi_t\right] dt  
\;+\; \si(t)\,\sqrt{( \xi_t)^+}\; d W_t  \;,\quad t\ge 0  
\eeqq
(note the particular form of the input term in comparison to (\ref{continuousequation}))  is given by 
$$
\left( J_{s,t}(\cdot,\cdot)\right)_{0\le s<t<\infty}    \quad,\quad J_{s,t}(y,\cdot) \;=\; H_{s,t}(y,\cdot) * I_{s,t}     
$$
where  $\left(I_{s,t}\right)_{s<t}$  is associated to $\al(\cdot)$ by lemma 3.6, and   $\left(H_{s,t}\right)_{s<t}$ is the semigroup of lemma 3.3. \\

{\bf Proof: } For $0\le s<t<\infty$ and $y\ge 0$, define $ J_{s,t}(y,\cdot):=H_{s,t}(y,\cdot) * I_{s,t}$. Combining (\ref{L6ii}) with the branching property according to lemma 2.2, we check that $\left( J_{s,t}(\cdot,\cdot) \right)_{0\le s<t<\infty}$ is a semigroup of transition probabilities on $\left( [0,\infty) , \calb([0,\infty)) \right)$. We have to verify that this semigroup corresponds to the process (\ref{continuousequationspecialrepresentation}). As in lemma 3.3, we give a heuristic argument. Since  $J_{s,t}(\cdot,\cdot)$ is a convolution,     
\beao
&&\left( J_{s,t}f - f \right)(y) \quad=\quad  \int J_{s,t}(y,dx)\, [f(x)-f(y)] \\  
&&=\quad   \int H_{s,t}(y,dv) \left\{ [f(v)-f(y)] \;+\; \int I_{s,t}(du)\, [f(v+u)-f(v)] \right\}  \;. 
\eeao
Consider the first expression on the right hand side: in case  $y=0$, this term is $0$ (since $H_{s,t}(0,\cdot)$ is the Dirac measure at $0$); in case  $y>0$, it behaves for $t\downarrow s$ as $\,(t-s)(\cala_s f)(y)\,$ in the notation of lemma 3.3~d). 
We exploit the moment structure  of $I_{s,t}$: from (\ref{L6i})+(\ref{L2i})+(\ref{L1aiii}) we have 
$$
E(I_{s,t}) \;=\; \int_s^t \frac{\si^2(v)}{2}\, \al(v)\, B(v,t)\,dv  \quad\sim\quad (t-s)\,\frac{\si^2(s)}{2}\, \al(s) 
\quad\quad\mbox{as}\quad t\downarrow s
$$
whereas  $\var\left(I_{s,t}\right)$ is proportional to $(t-s)^2$. In case $y=0$, the second expression above is 
$$
\int I_{s,t}(du)\, [f(y+u)-f(y)]  \quad\sim\quad (t-s)\,\frac{\si^2(s)}{2}\, \al(s)\, f'(y) \;+\; O\left( (t-s)^2 \right)  \quad,\quad t\downarrow s \;.  
$$
In case $y>0$, the second expression above is 
$$
\int I_{s,t}(du)\, [f(v+u)-f(v)]  \quad\sim\quad (t-s)\,\frac{\si^2(s)}{2}\, \al(s)\, f'(v) \;+\; O\left( (t-s)^2 \right)   \quad,\quad t\downarrow s
$$
where values $v$ selected by $H_{s,t}(y,dv)$ are concentrated on small neighbourhoods of their expectation  $y\,B(s,t)$, cf.\ (\ref{L3i}), hence close to $y$ as $t\downarrow s$. Hence in both cases $y=0$ (where by lemma~3.3 $(\cala_s f)(0)$ equals $0$) and $y>0$ we end up with  
$$
\lim_{t\downarrow s}  \frac{ J_{s,t}f - f }{t-s}(y)   \;\;=\;\;   \left\{ (\cala_s f)(y) + \frac{\si^2(s)}{2}\, \al(s)\, f'(y) \right\} 
\;,\quad y\ge 0
$$
which corresponds to equation (\ref{continuousequationspecialrepresentation}). \halmos \\

%
%
%

The next step is to specify the semigroup $\left( \wt J_{s,t}(\cdot,\cdot) \right)_{0\le s<t<\infty}$ associated with equation 
\beqq\label{onlyjumpinputequation}
d \xi_t \;=\; -\beta(t)\,\xi_{t^-}\, dt  \;+\; \int_{(0,\infty)} y\; \mu(dt,dy) 
\;+\; \si(t)\,\sqrt{(\xi_{t^-})^+}\; dW_t  \;,\quad t\ge 0     
\eeqq
where only the 'discrete' part of the input into equation $(\diamond)$ appears. We start with \\

{\bf 3.8 Lemma: } Let $\nu$ satisfy the restrictive condition (\ref{restrictivecondition}). Fix $0\le t_1<t_2\le t_3$ and write $\bar \mu (ds,dy)$ for the restriction of the random measure $\mu(ds,dy)$ to $(t_1,t_2]{\times}(0,\infty)$: then 
$$
\mbox{ $\bar \mu (ds,dy)$ is PRM on $(t_1,t_2]{\times}(0,\infty)$ with intensity $\bar \nu (ds,dy) := \wt a(s) 1_{ (t_1,t_2] }(s) ds\; \nu(dy)$ } \;. 
$$ 
Represent $\bar\mu$ in form $\sum_{i=1}^\infty \epsilon_{(T_i,Y_i)}$  and associate to  the point masses of $\bar\mu$ sitting at $(T_i,Y_i)$  processes $\left( \xi^{(T_i,Y_i)}_t \right)_{t\ge T_i}$ which are strong solutions to $(*,T_i,Y_i)$ according to lemma 2.1, driven by independent BM's and hence independent, conditionally on $\bar\mu$. Then thanks to  (\ref{restrictivecondition}), the sum 
\beqq\label{L7neuvar} 
\sum_{i=1}^\infty \,\xi^{(T_i,Y_i)}_{t_3}  \;=:\;   \int_{ (t_1,t_2]{\times}(0,\infty) } \xi_{t_3}^{(s,y)}\, \bar\mu (ds,dy) 
\eeqq
is finite a.s.\ and has the Laplace transform 
\beqq\label{L7neuvarbis} 
\la \;\;\lra\;\;  \exp\left\{ \;-\;  \int_{ (t_1,t_2]{\times}(0,\infty) } \left[ 1 -e^{ - y\, \Psi_{s,t_3}(\la)} \right]  \bar\nu(ds,dy)\, \right\} 
\eeqq
or with notation of (\ref{3ter})  
$$
\la \;\;\lra\;\;  \exp\left\{ \;-\;  \int_{t_1}^{t_2} \wt a(s)\, \wt \Psi_{s,t_3}(\la)\, ds \, \right\}  \quad,\quad \la\ge 0 \;. 
$$

\vskip0.5cm
{\bf Proof: } 1)~ Consider first some $a>0$ and some finite collection of points $(s_i,u_i)$ in $(t_1,t_2]{\times}[a,\infty)$, $1\le i\le \ell$. Associate to the points $(s_i,u_i)$ processes $\left( \xi^{(s_i,u_i)}_t \right)_{t\ge s_i}$ which are strong solutions to $(*,s_i,u_i)$, driven by independent BM's and hence independent. Then by lemma 3.3, the sum  
$$
\sum_{i=1}^\ell \,\xi^{(s_i,u_i)}_{t_3} 
$$
has  Laplace transform 
$$
\la \;\;\lra\;\;  \prod_{i=1}^\ell \;e^{\,-\, u_i\, \Psi_{s_i,t_3}(\la) } \;\;,\;\;\la\ge 0 \;. 
$$

2)~ Arrange the point masses $(T_i,Y_i)$ of $\bar\mu$ in decreasing order $\infty>Y_1>Y_2>\ldots$ of the second component. Let $\delta_n \downarrow 0$ denote the sequence selected in the proof of lemma 2.4, based on  (\ref{restrictivecondition}). Then for every $n$ fixed, we read step 1)  conditionally on $\bar\mu$: 
\beqq\label{hilf1}
E\left( e^{\,-\,\la\; \sum_{ i\in\bbn , Y_i\ge \delta_n } \xi^{(T_i,Y_i)}_{t_3} } \mid \calf_{\bar\mu} \right) \;\;=\;\;   \prod_{i\in\bbn , Y_i\ge \delta_n } \;e^{\,-\, Y_i\, \Psi_{T_i,t_3}(\la) } \;\;,\;\;\la\ge 0 \;. 
\eeqq
Under condition (\ref{restrictivecondition}) we have bounds (\ref{bound5}), and thus  as in the proof of lemma 2.4   
$$
\sum_{ i\in\bbn , Y_i\ge \delta_n } \xi^{(T_i,Y_i)}_{t_3}   \;\;\lra\quad  \sum_{i=1}^\infty \xi^{(T_i,Y_i)}_{t_3} \quad,\quad \nto 
$$
where the limit is a.s.\ finite. On the other hand, the collection of points $(Y_i)_{i\ge 1}$ being PRM on $(0,\infty)$ with intensity $[\int_{t_1}^{t_2}\wt a(s) ds] \cdot \nu(dy)$, $\;\sum_{i=1}^\infty Y_i$ is always finite by our permanent assumption $\int (y\wedge 1) \nu(dy) <\infty$. 
Moreover,   $s\to\Psi_{s,t_3}(\la)$ is bounded on $(t_1,t_2]$ for fixed $\la$ in virtue of lemma 3.1. Thus we have convergence on both sides   of (\ref{hilf1}) as $\nto$, and arrive at
\beqq\label{hilf2}
E\left( e^{\;-\, \la \; \sum_{i=1}^\infty \xi^{(T_i,Y_i)}_{t_3}} \mid \calf_{\bar\mu}\right) \;\;=\;\;  \prod_{i=1}^\infty \;e^{\,-\, Y_i\, \Psi_{T_i,t_3}(\la) }  \;\;,\;\;\la\ge 0 \;. 
\eeqq

3)~ Fix $\la\ge 0$ and define  $\,F : (t_1,t_2]{\times}(0,\infty) \to [0,\infty)\,$ by $\,F(s,y) = y\, \Psi_{s,t_3}(\la)\,$. Then we have 
$$
E\left( e^{ \,-\, \int F(s,y)\, \bar\mu(ds,dy) } \right)  \;=\;  
E\left( \exp\left\{ \,-\, \int_{ (t_1,t_2]{\times}(0,\infty) }  \left[ 1 - e^{\,-\, F(s,y)} \right] \bar\nu(ds,dy) \right\} \right) 
$$
by  [IW 89, p.\ 42--44]). Taking expectations in (\ref{hilf2}) we have proved (\ref{L7neuvarbis}).\halmos\\

{\bf 3.9 Lemma: } a)~ For all $0\le s<t<\infty$, there is a  law $\wt I_{s,t}$ on $[0,\infty)$ with Laplace transform 
\beqq\label{L8i}
\int_0^\infty  e^{-\la x}\, \wt I_{s,t}(dx)   \; \;=\; \;  \exp\left( -\int_s^t \wt a(v)\, \wt\Psi_{v,t}(\la)\, dv \right) \;\;,\;\; \la\ge 0
\eeqq
and the family $(\wt I_{s,t})_{0\le s<t<\infty}$ has the same skew convolution property 
\beqq\label{L8ii}
\wt I_{t_1,t_2}H_{t_2,t_3} * \wt I_{t_2,t_3} \;=\; \wt I_{t_1,t_3}  
\quad\mbox{for $0\le t_1<t_2<t_3<\infty$} \;. 
\eeqq
with respect to the semigroup $\left(H_{s,t}(\cdot)\right)_{0\le s<t<\infty}$ as the family $(I_{s,t})_{0\le s<t<\infty}$ of lemma 3.6. 

b)~ The semigroup  $( \wt J_{s,t}(\cdot,\cdot) )_{0\le s<t<\infty}$  of transition probabilities associated to the process (\ref{onlyjumpinputequation})   
$$
d \xi_t \;=\; -\beta(t)\,\xi_{t^-}\, dt  \;+\; \int_{(0,\infty)} y\; \mu(dt,dy) 
\;+\; \si(t)\,\sqrt{(\xi_{t^-})^+}\; dW_t  \;,\quad t\ge 0     
$$
is given by 
$$
\wt J_{s,t}(y,\cdot) \;=\; H_{s,t}(y,\cdot) * \wt I_{s,t}       \quad,\quad   0\le s<t<\infty \;,\; y\ge 0  \;. 
$$

\vskip0.8cm
{\bf Proof: } By definition of $\wt\Psi_{v,t}(\cdot)$ in (\ref{3ter}) and by 3.3 together with (\ref{L2iii}) , we have 
$$
\wt\Psi_{v,t_2}\circ \Psi_{t_2,t_3} \;=\; \wt\Psi_{t_2,t_3} \quad\mbox{for $t_1\le v\le t_2\le t_3$} \;, 
$$
thus the skew convolution property (\ref{L8ii}) has already been proved in lemma 3.8. This is a). 
 
As in lemma 2.2, we prepare independent BM's $W^{(i)}$ and  processes 

i)  $\;\xi^{(1)}=F^{(1)}(W^{(1)},s,y)$  strong solution to equation $(*,s,y)$ of lemma 2.1; by lemma 3.3:  
$$
\call\left(\, \xi^{(1)}_t \mid \xi^{(1)}_s=y \right)   \;=\;  H_{s,t}(y,\cdot) \;. 
$$

ii)  $\;\xi^{(2)}=F^{(2)}(W^{(2)},\mu, s,0)$  strong solution to equation (\ref{onlyjumpinputequation})   
starting in position $0$ at time $s$; combining lemmata 2.4 and 3.8 we do have 
$$
\call\left(\, \xi^{(2)}_t \mid \xi^{(2)}_s=0 \right)   \;=\; \wt I_{s,t} \;. 
$$
As strong solutions, both processes $\xi^{(1)} , \xi^{(2)}$ are independent, and  their sum $\xi^{(1)}+ \xi^{(2)}$ provides as in lemma 2.2 (this is again the 'branching property') a weak solution to equation  (\ref{onlyjumpinputequation}) starting in position $y$ at time $s$. An extension of the Yamada-Watanabe criterion, as in [H 08], gives pathwise uniqueness for equation (\ref{onlyjumpinputequation}). Thus assertion b) follows. \halmos\\

{\bf 3.10 Proof of theorem 1.2: }  Represent the unique strong solution $\xi$ of equation $(\diamond)$ with initial condition $y\ge 0$ at time $s$ as sum of two processes $\xi = \xi^{(1)}+\xi^{(2)}$ driven by independent Brownian motions $W^{(1)}$ and $W^{(2)}$, where $\xi^{(1)}=F(W^{(1)},s,y)$ is a strong solution to the CIR equation (\ref{continuousequation}) with time dependent coefficients 
$$
d \xi^{(1)}_t \;=\; \left[ a(t)-\beta(t)\, \xi^{(1)}_t\right] dt  
\;+\; \si(t)\,\sqrt{( \xi^{(1)}_t)^+}\; d W^{(1)}_t  \;,\quad  t\ge s \;,\; \xi^{(1)}_s=y 
$$
and $\xi^{(2)}=F(W^{(2)},\mu,s,0)$ a strong solution to  (\ref{onlyjumpinputequation}) starting in position $0$ at time $s$
$$
d \xi^{(2)}_t \;=\; -\beta(t)\,\xi^{(2)}_{t^-}\, dt  \;+\; \int_{(0,\infty)} y\; \mu(dt,dy) 
\;+\; \si(t)\,\sqrt{(\xi^{(2)}_{t^-})^+}\; dW^{(2)}_t  \;,\quad t\ge s  \;,\;  \xi^{(2)}_s=0  \;. 
$$
The branching property as in the last proof (or as in the proof 2.6 of theorem 1.1) combined with pathwise uniqueness allows to construct the  solution to $(\diamond)$ starting in position $y$ at time $s$ in form $\xi = \xi^{(1)}+\xi^{(2)}$, and independence  of 
$\xi^{(1)} , \xi^{(2)}$ together with lemmata 3.3 and 3.6+3.7 and 3.8+3.9 identifies the semigroup $\left(K_{s,t}(\cdot,\cdot)\right)_{0\le s<t<\infty}$ corresponding to the process $(\diamond)$ as  
\beqq\label{semigroupK}
K_{s,t}(y,\cdot)   \;\;=\;\;  H_{s,t}(y,\cdot) \;*\; \left( I_{s,t} * \wt I_{s,t} \right) \quad,\quad 0\le s<t<\infty \,,\, y\ge 0  
\eeqq
with Laplace transforms 
$$
\int_0^\infty e^{-\la x}\, K_{s,t}(y,dx)   \;\;=\;\;   
\exp\left( -\; y\,\Psi_{s,t}(\la) \;-\; \int_s^t  \left[\, a(v)\, \Psi_{v,t}(\la) + \wt a (v)\, \wt\Psi_{v,t}(\la) \,\right] dv \,\right) 
$$
as stated in theorem 1.2. The proof of theorem 1.2 is now complete.\halmos\\

\vskip1.0cm
{\large\bf References} \\
\small

[AN 72]\quad
Athreya, K,, Ney, P.: 
Branching processes. 
Springer 1972

\vskip0.2cm
[B 71]\quad
Barra, J.: 
Notions fondamentales de statistique mathematique. 
Dunod 1971. 

\vskip0.2cm
[BGT 89]\quad
Bingham, N., Goldie, C., Teugels, J.: 
Regular variation. 
Cambridge 1989. 

\vskip0.2cm
[B 81]\quad
Br\'emaud, P.: 
Point processes and queues. 
Springer 1981

\vskip0.2cm
[BH 06]\quad 
Brodda, K., H\"opfner, R.: 
A stochastic model and a functional central limit theorem for information processing in large systems of neurons. 
J. Math. Biol. {\bf 52}, 439--457 (2006). 

\vskip0.2cm
[CG 06]\quad
Chaleyat-Maurel, M., Genon-Catalot, V.: 
Computable infinite dimensional filters with applications to discretized diffusion processes. \\
Stoch. Process. Applications {\bf 116}, 1447--1467 (2006).

\vskip0.2cm
[CIR 85]\quad 
Cox, J., Ingersoll, J., Ross, S.: 
A theory of the term structure of interest rates. \\
Econometrica {\bf 53}, 385--407 (1985). 


\vskip0.2cm
[F 71]\quad
Feller, W.: 
An introduction to probability theory and its applications, Vol. II. \\
Wiley 1971. 

\vskip0.2cm
[GY 03]\quad
G\"oing-Jaeschke, A., Yor, M.: 
A survey and some generalizations of Bessel processes. \\
Bernoulli {\bf 9}, 313--350 (2003). 

\vskip0.2cm
[H 06]\quad
Hammer, M.: 
Parametersch\"atzung in zeitdiskreten ergodischen Markov-Prozessen am Beispiel des Cox-Ingersoll-Ross Modells. 
Diplomarbeit, Institute of Mathematics, University of Mainz, 2006 (see 
{\tt http://archimed.uni-mainz.de/opusubm/archimed-home.html}). 

\vskip0.2cm
[H 07]\quad 
H\"opfner, R.: 
On a set of data for the membrane potential in a neuron. \\
Math.\ Biosci.\ {\bf 207}, 275--301 (2007)

\vskip0.2cm
[H 08]\quad 
H\"opfner, R.: 
An extension of the Yamada-Watanabe condition for pathwise unicity to stochastic differential equations with jumps. \\
Preprint 2008 (see {\tt www.mathematik.uni-mainz.de/$\sim$hoepfner}), submitted. 

\vskip0.2cm
[HJ 94]\quad 
H\"opfner, R., Jacod, J.: 
Some remarks on the joint estimation of the index and the scale parameter for stable processes.
In: Mandl, P., Huskova, M. (Eds.): Asymptotic Statistics.
Proceedings of the Fifth Prague Symposium 1993, pp. 273-284, 
Physica Verlag 1994. 

\vskip0.2cm
[HK 09]\quad 
H\"opfner, R., Kutoyants, Yu.: Estimating discontinuous periodic signals in a time inhomogeneous diffusion. 
Preprint 2009 (see {\tt www.mathematik.uni-mainz.de/$\sim$hoepfner}), submitted. 

\vskip0.2cm
[IW 89]\quad 
Ikeda, N., Watanabe, S.: 
Stochastic differential equations and diffusion processes. \\
2nd ed., North-Holland, 1989. 

\vskip0.2cm
[J 09]\quad 
Jahn, P.: Statistical Problems related to excitation threshold and reset value of membrane potentials. PhD thesis, Institute of Mathematics, University of Mainz, 2009. 

\vskip0.2cm
[JS 87]\quad 
Jacod, J., Shiryaev, A.: 
Limit theorems for stochastic processes. 
Springer 1987. 

\vskip0.2cm
[KS 91]\quad
Karatzas, I., Shreve, S.: 
Brownian motion and stochastic calculus. \\
2nd ed. Springer 1991. 

\vskip0.2cm 
[LL 87]\quad
L\'ansk\'y, P., L\'ansk\'a, V.: Diffusion approximation of the neuronal model with synaptic reversal potentials. 
Biol.\ Cybern.\ {\bf 56}, 19--26 (1987). 


\vskip0.2cm
[L 03]\quad 
Li, Z.: 
Skew convolution semigroups and related immigration processes. \\
Theory Probab. Appl. {\bf 46}, 274--296 (2003). 


\vskip0.2cm
[OR 97]\quad 
Overbeck, L., Ryd\'en, T.: 
Estimation in the Cox-Ingersoll-Ross model. \\
Econometric Theory {\bf 13}, 430--461 (1997). 

\vskip0.2cm
[P 05]\quad 
Protter, P.: 
Stochastic integration and differential equations. \\
Springer 1990, 2nd ed.\ 2005. 

\vskip0.2cm
[YW 71]\quad 
Yamada, T., Watanabe, S.: 
On the uniqueness of solutions of stochastic differential equations. 
J.\ Math.\ Kyoto Univ.\ {\bf 11}, 155-167 (1971).

\normalsize
\vskip1.5cm
~\hfill {\bf 09.03.2009}

\small
\vskip0.5cm
Reinhard H\"opfner\\
Institut f\"ur Mathematik, Universit\"at Mainz, D--55099 Mainz\\ 
{\tt hoepfner@mathematik.uni-mainz.de}\\
{\tt http://www.mathematik.uni-mainz.de/$\sim$hoepfner}

\end{document}